\documentclass[12pt]{article}
\usepackage{amsmath, amssymb, amsthm}
\usepackage{geometry}
\usepackage{hyperref}
\usepackage{enumitem}
\usepackage{mathrsfs}
\usepackage{cite}
\newtheorem{theorem}{Theorem}[section]
\newtheorem{lemma}[theorem]{Lemma}
\newtheorem{corollary}[theorem]{Corollary}
\newtheorem{proposition}[theorem]{Proposition}

\newtheorem*{example*}{Example}
\newtheorem*{remark*}{Remark}

\title{The Basic Component of the Mean Curvature of Riemannian Foliations }

\author{Jes\'us A. \'Alvarez L\'opez\\
Universidade de Santiago de Compostela\\
Departamento de Xeometr\'ia e Topolox\'ia\\
Colexio Universitario de Lugo\\
27071 Lugo, Spain}
\date{}

\begin{document}
\maketitle

\begin{abstract}
For a Riemannian foliation $\mathcal{F}$ on a compact manifold $M$ with a bundle-like metric, the de Rham complex of $M$ is $\mathcal{C}^{\infty}$-splitted as the direct sum of the basic complex and its orthogonal complement. Then the basic component $\kappa_{b}$ of the mean curvature form of $\mathcal{F}$ is closed and defines a class $\xi(\mathcal{F})$ in the basic cohomology that is invariant under any change of the bundle-like metric. Moreover, any element in $\xi(\mathcal{F})$ can be realized as the basic component of the mean curvature of some bundle-like metric.

It is also proved that $\xi(\mathcal{F})$ vanishes iff there exists some bundle-like metric on $M$ for which the leaves are minimal submanifolds. As a consequence, this tautness property is verified in any of the following cases: (a) when the Ricci curvature of the transverse Riemannian structure is positive, or (b) when $\mathcal{F}$ is of codimension one. In particular, a compact manifold with a Riemannian foliation of codimension one has infinite fundamental group.

\end{abstract}

\section*{Introduction}
The basic complex of a Riemannian foliation $\mathcal{F}$ on a compact manifold $M$ has been extensively studied since it was defined by Reinhart in [24]. One of the main interests was to compare its properties with the properties of the de Rham complex of a compact Riemannian manifold.

There have been two main ways of doing this study. The first one depends strongly on the structure theorems of P. Molino ([21]). It consists of reducing problems of Riemannian foliations to problems of transversally parallelizable foliations, or Lie foliations with dense leaves, where they have easier solutions. It is important that this method depends only on the transverse structure of $\mathcal{F}$. In fact, every result can be equally proved for the invariant de Rham complex of a complete pseudogroup of local isometries ([13]) with a compact space of orbit closures (using the results of [27]).

This method has been used to study for instance the finite dimension and duality of the basic cohomology $H_{b}(\mathcal{F})$ ([9], [15], [28]), and also the Hodge decomposition of the basic complex $A_{b}(\mathcal{F})$ with respect to some scalar product defined by using Molino's structure theorems ([7], [8]). Especially interesting is that, even when $\mathcal{F}$ is transversally oriented, duality not always holds in $H_{b}(\mathcal{F})$, which is different from the situation for the de Rham cohomology of a compact manifold. With this method it was obtained that $H_{b}(\mathcal{F})$ has duality iff $H_{b}^{q}(\mathcal{F}) \neq 0$ ($q=\operatorname{codim}(\mathcal{F})$) ([8]). On the other hand, also with this method, X.~Masa has proved the following tautness theorem ([19]): when $\mathcal{F}$ and $M$ are oriented, $H_{b}^{q}(\mathcal{F}) \neq 0$ iff $\mathcal{F}$ is taut, i.e. iff there exists a Riemannian metric on $M$ for which all the leaves are minimal submanifolds. (This theorem was also proved in [10], [12], [13], [22], and [2] for special cases.) So the mean curvature of $\mathcal{F}$ has a certain influence on the duality of $H_{b}(\mathcal{F})$.

The second way of studying the basic cohomology uses strongly the mean curvature of $\mathcal{F}$. For instance a basic Hodge decomposition was obtained with respect to the restriction to $A_{b}(\mathcal{F})$ of the scalar product of the de Rham complex $A(M)$ of $M$ ([18]). Hence this second decomposition gives a more direct relation between the geometry of $(M, \mathcal{F})$ and the basic cohomology, but it requires the additional hypothesis for the mean curvature form to be basic. The proof gives an explicit description of the basic Laplacian using the mean curvature, so that it can be extended to an elliptic operator on $A(M)$. This is the reason of the crucial role that the mean curvature plays in this method, and it shows how the peculiar features of the basic cohomology depend strongly on the mean curvature ([17], [30]). For example, when the mean curvature is basic, the above tautness theorem follows directly from this basic Hodge decomposition.

In the present paper we basically show how to extend the second method above to general Riemannian foliations on compact manifolds. The main idea is to obtain a $\mathcal{C}^{\infty}$-splitting of $A(M)$ so that the mean curvature form can be orthogonally projected to the basic complex. This basic component of the mean curvature will be the key point of our study.

The indicated splitting of $A(M)$ is really a partial leafwise Hodge decomposition. In [3] a leafwise Hodge decomposition for Riemannian foliations on compact manifolds was proven. It works for the space of $L_{2}$-differential forms that are smooth along the leaves. In that paper it was also conjectured that this decomposition still holds for smooth forms with the $\mathcal{C}^{\infty}$-topology. The first result of our paper is that this conjecture is true at least for horizontal forms; i.e. $A(M)$ can be $\mathcal{C}^{\infty}$-decomposed as the direct sum of $A_{b}(\mathcal{F})$ and its orthogonal complement, $A_{b}(\mathcal{F})^{\perp}$. A description of $A_{b}(\mathcal{F})^{\perp}$ using the codifferential map along the leaves is also given.

This decomposition of $A(M)$ implies that the restriction of the exterior derivative to $A_{b}(\mathcal{F})$ has an adjoint operator, $\delta_{b}$, given by orthogonal projection of the codifferential map of $M$. We describe $\delta_{b}$ using the $A_{b}(\mathcal{F})$-component of the mean curvature form, denoted by $\kappa_{b}$. It turns out that $\kappa_{b}$ is closed, defining a class $\left[\kappa_{b}\right]$ in $H_{b}^{1}(\mathcal{F})$.

We further prove that $\left[\kappa_{b}\right]$ is invariant under any change of the bundle-like metric defining any transverse Riemannian structure. So $\left[\kappa_{b}\right]$ is an invariant of the foliation, and will be denoted by $\xi(\mathcal{F})$. Moreover, any element of $\xi(\mathcal{F})$ can be realized as the basic component of the mean curvature of some bundle-like metric defining a given transverse Riemannian structure. The proof of these results is reduced to two cases:\\
(i) Modifying only the orthogonal complement of the leaves, and\\
(ii) Modifying only the metric along the leaves.

In both cases the results follow by direct computation using Rummler's mean curvature formula ([26]). Furthermore $\kappa_{b}$ is also invariant under type (i) changes of the metric.

The form $\kappa_{b}$ can be used to prove a Hodge decomposition of $A_{b}(\mathcal{F})$ with respect to the restriction of the scalar product of $A(M)$. In fact, the methods of [18] and [23] can be generalized without the hypothesis on the mean curvature to be basic. As a consequence, when $\mathcal{F}$ is transversally oriented, $H_{b}(\mathcal{F})$ has duality iff $\xi(\mathcal{F})=0$. This duality characterization and the tautness theorem of X.~Masa imply that $\xi(\mathcal{F})$ vanishes iff $\mathcal{F}$ is taut. So $\operatorname{Ric}_{T}>0$ implies that $\mathcal{F}$ is taut by a result of Hebda ([14]), where $\operatorname{Ric}_{T}$ is the Ricci curvature of the corresponding invariant metric, $g_{T}$, on the holonomy pseudogroup of $\mathcal{F}$. Another consequence is that every Riemannian foliation of codimension one on a compact manifold is taut, obtaining that compact manifolds with Riemannian foliations of codimension one have infinite fundamental group.

The class $\xi(\mathcal{F})$ was also defined in [6] in a very different way, and without relating it to the mean curvature.

After finishing this paper, the author has been informed that F.~Kamber has proved the following: for any Riemannian foliation on a compact manifold there exists some bundle-like metric with basic mean curvature. Thus our study for that special metric reduces to the case of [16], [17], [18] and [23]. Nevertheless, this paper hopefully sheds some more light on the interplay of geometrical and cohomological properties of arbitrary bundle-like metrics.

The author would like to thank X.~Masa and Ph.~Tondeur for helpful conversations.

\section{Preliminaries}
Let $M$ be a smooth manifold of dimension $n$, and let $(A(M), d)$ be its de Rham complex. Let $\mathcal{F}$ be a smooth foliation of dimension $p$ and codimension $q$ on $M$. Then we have the vector subbundle $T \mathcal{F} \subset T M$ of vectors tangent to the leaves, and let $\mathcal{X}(\mathcal{F})=\Gamma(T \mathcal{F})$. Let also $\mathcal{X}(M, \mathcal{F}) \subset \mathcal{X}(M)$ be the subspace of infinitesimal transformations of $(M, \mathcal{F})$.

The complex $(A_{b}(\mathcal{F}), d_{b})$ of basic differential forms is defined as
\begin{equation*}
A_{b}(\mathcal{F})=\left\{\alpha \in A(M) \mid \iota_{X} \alpha=\theta_{X} \alpha=0 \quad \text { for } \quad X \in \mathcal{X}(\mathcal{F})\right\}, \tag{1.1}
\end{equation*}
with $d_{b}=d_{\mid A_{b}(\mathcal{F})}([24])$.

For a given Riemannian metric on $M$, $g=(\cdot, \cdot)$, we have the orthogonal decomposition
\begin{equation*}
T M=Q \oplus T \mathcal{F}, \tag{1.2}
\end{equation*}
where $Q=T \mathcal{F}^{\perp}$. Then we obtain the associated bigrading of $A(M)$ given by
\begin{equation*}
A^{u, v}(M)=\Gamma\left(\Lambda^{u} Q^{*} \otimes \Lambda^{v} T^{*} \mathcal{F}\right)=\Gamma\left(\Lambda^{u} Q^{*}\right) \otimes_{\mathcal{C}^{\infty}(M)} \Gamma\left(\Lambda^{v} T^{*} \mathcal{F}\right) \tag{1.3}
\end{equation*}
The exterior derivative decomposes as $d=d_{0,1}+d_{1,0}+d_{2,-1}$, where each $d_{i, j}$ is bihomogeneous of bidegree $(i, j)$. These bihomogeneous components verify well-known properties ([1]). For example
\begin{equation*}
A_{b}(\mathcal{F})=A^{., 0}(M) \cap \operatorname{Ker}\left(d_{0,1}\right). \tag{1.4}
\end{equation*}
Clearly we also have
\begin{equation*}
d_{0,1} \equiv d_{\mathcal{F}} \quad \text { on } \quad A^{0, \cdot}(M) \equiv \Gamma\left(\Lambda T^{*} \mathcal{F}\right),\tag{1.5}
\end{equation*}
where $d_{\mathcal{F}}$ is the exterior derivative along the leaves.

The formal adjoint of $d$ also decomposes as $\delta=\delta_{0,-1}+\delta_{-1,0}+\delta_{-2,1}$ analogously, where each $\delta_{i, j}$ is the formal adjoint of $d_{-i,-j}$.

\begin{lemma}
$d_{2,-1}$ and $\delta_{-2,1}$ are differential operators of order zero.
\end{lemma}

\begin{proof} It is enough to check the case of $d_{2,-1}$, which follows from the formula for the exterior derivative of a product since $d_{2,-1}\left(\mathcal{C}^{\infty}(M)\right)=0$.
\end{proof}

Suppose now that $\mathcal{F}$ is oriented and transversally oriented, and let $*_{Q}$ and $*_{\mathcal{F}}$ be the corresponding star operators on $\Gamma\left(\Lambda Q^{*}\right)$ and $\Gamma\left(\Lambda T^{*} \mathcal{F}\right)$, respectively. Choose the orientation of $M$ such that $(X_{1}, \ldots, X_{q}, Y_{1}, \ldots, Y_{p})$ is an oriented frame of $M$ if $(X_{1}, \ldots, X_{q})$ and $(Y_{1}, \ldots, Y_{p})$ are oriented frames of $Q$ and $T \mathcal{F}$ respectively. With this choice of orientations, it can be proved as in [3] that
\begin{equation*}
*=(-1)^{(q-u) v} *_{Q} \otimes *_\mathcal{F} \quad \text { on } \quad A^{u, v}(M) \tag{1.6}
\end{equation*}
according to the tensor decomposition (1.3). (In [3] the signs are different because of the different choice of orientations.) Extend the operators $*_Q$ and $*_{\mathcal{F}}$ to $A(M)$ by defining $*_{Q}=*_{Q} \otimes\operatorname{id}$ and $*_{\mathcal{F}}=\operatorname{id}\otimes *_{\mathcal{F}}$.

From now on assume that $\mathcal{F}$ is Riemannian and the metric bundle-like ([24]). Let $\nabla$ be the Levi-Civita connection of $M$, and let $\stackrel{\circ}{\nabla}$ be the connection on $T M$ defined as the direct sum of the connections that $\nabla$ induces on $Q$ and $T \mathcal{F}$ ([3]). $\stackrel{\circ}{\nabla}$ induces a connection on $\Lambda T^{*} M$, also denoted by $\stackrel{\circ}{\nabla}$, such that each $\stackrel{\circ}{\nabla}_{X}(X \in$ $\mathcal{X}(M)$ ) is the bihomogeneous component of bidegree $(0,0)$ of $\nabla_{X}$. Let $E_{1}, \ldots, E_{p}$ be a local orthonormal frame of $T \mathcal{F}$, and $\beta_{1}, \ldots, \beta_{p}$ the dual coframe. Then, for any $\omega \in A^{0, \cdot( } M$ ), we have locally ([3])
\begin{align*}
& d_{0,1}(\omega)=\sum_{i=1}^{p} \beta_{i} \wedge \stackrel{\circ}{\nabla}_{E_{i}} \omega,  \tag{1.7}\\
& \delta_{0,-1}(\omega)=-\sum_{i=1}^{p} \iota_{E_{i}} \stackrel{\circ}{\nabla}_{E_{i}} \omega . \tag{1.8}
\end{align*}
Therefore we also have
\begin{equation*}
\delta_{0,-1} \equiv \delta_{\mathcal{F}} \quad \text { on } \quad A^{0, \cdot}(M) \equiv \Gamma\left(\Lambda T^{*} \mathcal{F}\right), \tag{1.9}
\end{equation*}
where $\delta_{\mathcal{F}}$ is the leafwise codifferential operator.

On a distinguished chart $U$ we can find an orthonormal frame $\alpha_{1}, \ldots, \alpha_{q}$ of $Q_{| U}^{*}$ such that each $\alpha_{i}$ is basic. Then the wedge product gives the identity.
\begin{equation*}
A^{u, v}(U) \equiv \Lambda^{u}\Big(\sum_{i=1}^{q} \mathbb{R} \cdot \alpha_{i}\Big) \otimes A^{0, v}(U). \tag{1.10}
\end{equation*}
As in [3], it can be proved that
\begin{align*}
& d_{0,1}=(-1)^{u} \text { id } \otimes d_{0,1},  \tag{1.11}\\
& \delta_{0,-1}=(-1)^{u} \text { id } \otimes \delta_{0,-1}, \tag{1.12}
\end{align*}
with respect to (1.10). Hence, by (1.5), (1.9), (1.11) and (1.12),
\begin{equation*}
\delta_{0,-1}=(-1)^{p v+p+1} *_\mathcal{F} d_{0,1} *_\mathcal{F} \quad \text { on } \quad A^{\cdot, v}(M). \tag{1.13}
\end{equation*}

From now on assume also that $M$ is compact. Consider on $A(M)$ the $\mathcal{C}^{\infty}$-topology, and on the space of currents, $A(M)^{\prime}$, consider the topology of uniform convergence on bounded subsets of $A(M)$. Let $\mathcal{D}: A(M) \longrightarrow A(M)^{\prime}$ be the de Rham duality map defined by
\begin{equation*}
\mathcal{D}(\alpha)(\beta)=\int_{M} \alpha \wedge \beta . \tag{1.14}
\end{equation*}

For each integer $v$ let
\begin{align*}
& \mathcal{Z}_{v}=\left(\operatorname{Ker}\left(d_{0,1}\right) \cap A^{,, v}(M)\right) \oplus \bigoplus_{j<v} A^{, j}(M),  \tag{1.15}\\
& \mathcal{B}_{v}=d_{0,1}\left(A^{, v-1}(M)\right) \oplus \bigoplus_{j<v} A^{, j}(M) . \tag{1.16}
\end{align*}

We have
\begin{equation*}
\mathcal{D}\left(\mathcal{Z}_{v}\right)=\left\{T \in \mathcal{D}(A(M)) \mid T\left(\mathcal{B}_{p-v}\right)=0\right\} \equiv\left(A(M) / \overline{\mathcal{B}_{p-v}}\right)^{\prime} \cap \mathcal{D}(A(M)), \tag{1.17}
\end{equation*}
which is dense in $\left(A(M) / \overline{\mathcal{B}_{p-v}}\right)^{\prime}$ by Proposition 5.1 of [1].

\begin{lemma} We have\\
{\rm i)} $\mathcal{B}_{p-v}^{\perp}=* \mathcal{Z}_{v}$;\\
{\rm ii)} $\mathcal{B}_{p-v}^{\perp \perp}=\overline{\mathcal{B}_{p-v}}$.
\end{lemma}

\begin{proof} i) follows easily from (1.14) and (1.17). On the other hand we have clearly
$$
\overline{\mathcal{B}_{p-v}} \subset \mathcal{B}_{p-v}{ }^{\perp \perp} .
$$
Now take $\alpha \in \mathcal{B}_{p-v}{ }^{\perp \perp}=\left(* \mathcal{Z}_{v}\right)^{\perp}$. Then for all $\beta \in \mathcal{Z}_{v}$ we have
$$
0=\langle\alpha, * \beta\rangle= \pm \mathcal{D}(\beta)(\alpha)
$$
So $T(\alpha)=0$ for all $T \in\left(A(M) / \overline{\mathcal{B}_{p-v}}\right)^{\prime}$, which implies that $\alpha \in \overline{\mathcal{B}_{p-v}}$ by the Hahn-Banach theorem.
\end{proof}

It is easy to check that
\begin{equation*}
*_Q A_{b}(\mathcal{F})=A_{b}(\mathcal{F}). \tag{1.18}
\end{equation*}
So
\begin{equation*}
*_{Q}\left(A_{b}(\mathcal{F})^{\perp}\right)=A_{b}(\mathcal{F})^{\perp}, \tag{1.19}
\end{equation*}
and by (1.6) and (1.11) we obtain
\begin{equation*}
* \mathcal{Z}_{v}=*_\mathcal{F} \mathcal{Z}_{v}. \tag{1.20}
\end{equation*}

\section{The Orthogonal Complement of the Basic Complex}
In this section we will prove the following result.

\begin{theorem} If $\mathcal{F}$ is a Riemannian foliation on a compact Riemannian manifold $M$, then
\begin{equation*}
A(M)=A_{b}(\mathcal{F}) \oplus A_{b}(\mathcal{F})^{\perp} \tag{2.1}
\end{equation*}
with the $\mathcal{C}^{\infty}$-Fréchet topology. Moreover, if the metric is bundle-like then
\begin{equation*}
A_{b}(\mathcal{F})^{\perp}=\overline{\delta_{0,-1}(A^{\cdot, 1}(M))} \oplus \bigoplus_{v>0} A^{\cdot, v}(M). \tag{2.2}
\end{equation*}
\end{theorem}

To prove (2.1) it is enough to consider any fixed bundle-like metric on $M$.

Assume firstly that $\mathcal{F}$ is transversally parallelizable (T.P.) ([21]). Let $X_{1}, \ldots, X_{q} \in$ $\mathcal{X}(M, \mathcal{F}) \cap \Gamma(Q)$ define an orthonormal transverse parallelism, and let $V \subset \mathcal{X}(M, \mathcal{F})$ be the subspace generated by the $X_{i}$ 's.

If $\pi_{b}: M \longrightarrow W$ is the basic fibering of $\mathcal{F}$ ([21]), then we have
\begin{align*}
& A_{b}^{0}(\mathcal{F})=\pi_{b}^{*}\left(\mathcal{C}^{\infty}(W)\right) , \tag{2.3}\\
& A^{u, v}(M) \equiv \Lambda^{u} V^{*} \otimes A^{0, v}(M). \tag{2.4}
\end{align*}
Thus it is enough to prove
\begin{equation*}
\mathcal{C}^{\infty}(M)=\pi_{b}^{*}\left(\mathcal{C}^{\infty}(W)\right) \oplus\left(\pi_{b}^{*}\left(\mathcal{C}^{\infty}(W)\right)\right)^{\perp}. \tag{2.5}
\end{equation*}
We can suppose that $M$ and $W$ are oriented, otherwise we can reduce the problem to the corresponding orientation coverings. Then let $\bar{\chi}$ be the characteristic form ([26]) of the foliation $\overline{\mathcal{F}}$ defined by $\pi_{b}$, whose leaves are the closures of the leaves of $\mathcal{F}$. We can define the continuous operator $\rho_{\mathcal{F}}: \mathcal{C}^{\infty}(M) \longrightarrow \mathcal{C}^{\infty}(W)$ by
\begin{equation*}
\rho_{\mathcal{F}}(f)=\frac{1}{h} \int_{\pi_{b}} f \bar{\chi}, \tag{2.6}
\end{equation*}
where $h=\int_{\pi_{b}} \bar{\chi}>0$. Clearly $\rho_{\mathcal{F}} \pi_{b}*=\mathrm{id}$, so
\begin{equation*}
\mathcal{C}^{\infty}(M)=\pi_{b}^{*}\left(\mathcal{C}^{\infty}(W)\right) \oplus \operatorname{Ker}\left(\rho_{\mathcal{F}}\right) .\tag{2.7}
\end{equation*}
Moreover, for $f \in \operatorname{Ker}\left(\rho_{\mathcal{F}}\right)$ and $g \in \mathcal{C}^{\infty}(W)$ we have
$$
\left\langle\pi_{b}^{*}(g), f\right\rangle=\int_{M} \pi_{b}^{*}(g) f \operatorname{Vol}_{M}=\int_{W} g h \rho_{\mathcal{F}}(f) \operatorname{Vol}_{W}=0,
$$
from which we have
\begin{equation*}
\operatorname{Ker}\left(\rho_{\mathcal{F}}\right)=\left(\pi_{b}^{*}\left(\mathcal{C}^{\infty}(W)\right)\right)^{\perp}. \tag{2.8}
\end{equation*}
Therefore, (2.5) follows from (2.7) and (2.8). Thus we obtain (2.1) when $\mathcal{F}$ is T.P.

In the general case we consider the $O(q)$-principal bundle of transverse orthonormal frames, $\pi: \widehat{M} \rightarrow M$, with the transverse Levi-Civita connection ([21]). Let $(o(q), \iota, \theta, A(\widehat{M}), d)$ be the associated action with the corresponding algebraic connection ([11]). On $\widehat{M}$ we consider the $O(q)$-invariant metric defined as the orthogonal sum of the lifting of the metric of $M$ on the horizontal bundle, and any metric on the vertical bundle for which some basis of fundamental fields is orthonormal. We can suppose that the fibres of $\pi$ have volume one. This is a bundle-like metric for the horizontal lifting, $\widehat{\mathcal{F}}$, of $\mathcal{F}$ to $\widehat{M}$.

Let $\chi_{\pi}$ be the characteristic form of the foliation defined by $\pi$, and let $\rho_{\pi}: A(\widehat{M}) \rightarrow$ $A(M)$ be the continuous operator defined by
\begin{equation*}
\rho_{\pi}(\alpha)=\int_{\pi} \alpha \wedge \chi_{\pi}. \tag{2.9}
\end{equation*}
Then it is easy to check that $\rho_{\pi} \pi^{*}=\mathrm{id}$, and
\begin{equation*}
\left\langle\rho_{\pi}(\alpha), \beta\right\rangle=\left\langle\alpha, \pi^{*} \beta\right\rangle \tag{2.10}
\end{equation*}
for all $\alpha \in A(\widehat{M})$ and $\beta \in A(M)$.

$\widehat{\mathcal{F}}$ has an orthonormal transverse parallelism defined by $X_{1}, \ldots, X_{q}, Y_{1}, \ldots, Y_{q_{0}} \in$ $\mathcal{X}(\widehat{M}, \widehat{\mathcal{F}})\left(q_{0}=q(q-1) / 2\right)$, where the $X_{i}$'s are horizontal and the $Y_{j}$'s are vertical. Let $V_{1} \subset \mathcal{X}(\widehat{M}, \widehat{\mathcal{F}})$ be the subspace generated by the $X_{i}$'s. We have
\begin{equation*}
A^{\cdot, 0}(\widehat{M})_{\iota=0}=\Lambda V_{1}^{*} \otimes \mathcal{C}^{\infty}(\widehat{M}). \tag{2.11}
\end{equation*}
Then, by tensoring $\rho_{\widehat{\mathcal{F}}}: \mathcal{C}^{\infty}(\widehat{M}) \rightarrow A_{b}^{0}(\widehat{\mathcal{F}})$ with the identity on $\Lambda V_{1}^{*}$ we get a retraction
$$
\rho_{\widehat{\mathcal{F}}}: A^{\cdot, 0}(\widehat{M})_{\iota=0} \rightarrow A_{b}(\widehat{\mathcal{F}})_{\iota=0}.
$$
So $\rho_{\mathcal{F}}=\rho_{\pi} \rho_{\widehat{\mathcal{F}}}$ defines a retraction of $A^{\cdot, 0}(\widehat{M})_{\iota=0, \theta=0} \equiv A^{\cdot, 0}(M)$ onto $A_{b}(\widehat{\mathcal{F}})_{\iota=0, \theta=0} \equiv$ $A_{b}(\mathcal{F})$, obtaining
\begin{equation*}
A^{0, \cdot}(M)=A_{b}(\mathcal{F}) \oplus \operatorname{Ker}\left(\rho_{\mathcal{F}}\right). \tag{2.12}
\end{equation*}
But (2.8) and (2.10) imply
\begin{equation*}
\operatorname{Ker}\left(\rho_{\mathcal{F}}\right)=A_{b}(\mathcal{F})^{\perp} \cap A^{\cdot, 0}(M). \tag{2.13}
\end{equation*}
So we obtain (2.1) from (2.12) and (2.13).

By Lemma 1.2 and (1.20) we have
\begin{equation*}
\left(*_{\mathcal{F}} A_{b}(\mathcal{F})\right)^{\perp}=\overline{\mathcal{B}_{p}}. \tag{2.14}
\end{equation*}
Then (2.2) follows from (1.13) and (2.14), which finishes the proof of Theorem 2.1.

For any $\alpha \in A(M)$ we will denote by $\alpha_{b}$ its orthogonal projection to $A_{b}(\mathcal{F})$, and by $\alpha_{o}$ its orthogonal projection to $A_{b}(\mathcal{F})^{\perp}$; i.e. $\alpha_{b}=\rho_{\mathcal{F}}(\alpha)$ and $\alpha_{o}=\alpha-\rho_{\mathcal{F}}(\alpha)$. Then the following result follows from the definition of $\rho_{\mathcal{F}}$.

\begin{proposition}$\ $\newline
{\rm i)} For $f \in \mathcal{C}^{\infty}(M), f>0$ implies $f_{b}>0$.\\
{\rm ii)} If $\alpha \in A(M)$ and $\beta \in A_{b}(\mathcal{F})$, then $(\alpha \wedge \beta)_{b}=\alpha_{b} \wedge \beta$.
\end{proposition}

\section{The Basic Component of the Mean Curvature}
In the above situation suppose that $\mathcal{F}$ is oriented and transversally oriented. Then, if $\nu \in A_{b}^{q}(\mathcal{F})$ is the transverse volume form, the characteristic form of $\mathcal{F}$ is $\chi=* \nu$. We also have Rummler's formula ([26])
\begin{equation*}
d_{1,0}(\chi)=-\kappa \wedge \chi, \tag{3.1}
\end{equation*}
where $\kappa$ is the mean curvature form of $\mathcal{F}$.

For $\alpha \in A^{r}(M)$ let $\alpha \neg$ denote the operator on $A(M)$ defined by (see [25])
\begin{equation*}
\alpha \neg \beta=\varepsilon(r, s) *(\alpha \wedge * \beta) \quad \text { if } \quad \beta \in A^{s}(M), \tag{3.2}
\end{equation*}
where
\[
\varepsilon(r, s)= \begin{cases}(-1)^{n s+n+1} & \text { if } r \equiv 1(\bmod 4)  \tag{3.3}\\ (-1)^{n s+s} & \text { if } r \equiv 2(\bmod 4) \\ (-1)^{n s+n} & \text { if } r \equiv 3(\bmod 4) \\ (-1)^{n s+s+1} & \text { if } r \equiv 0(\bmod 4)\end{cases}.
\]
If $X \in \Gamma\left(\Lambda^{r} T M\right)$ is $(\cdot, \cdot)$-dual to $\alpha$ (i.e. $\alpha=(X, \cdot)$ ), then it is easy to check by induction on $r$ that
\begin{equation*}
\alpha \neg \beta=-\iota_{X} \beta. \tag{3.4}
\end{equation*}
If $\alpha \in A^{r, 0}(M)$ and $\beta \in A^{s, 0}(M)$ we also have
\begin{equation*}
\alpha \neg \beta=\varepsilon(r, s) * Q\left(\alpha \wedge *_{Q} \beta\right) . \tag{3.5}
\end{equation*}
Thus, by (1.18),
\begin{equation*}
\alpha \neg A_{b}(\mathcal{F}) \subseteq A_{b}(\mathcal{F}) \quad \text { if } \quad \alpha \in A_{b}(\mathcal{F}). \tag{3.6}
\end{equation*}
Moreover, the following properties have straightforward proofs:
\begin{align*}
& (\alpha \wedge \beta, \gamma)=-(-1)^{r(r-1) / 2}(\beta, \alpha \neg \gamma) \quad \text { for } \quad \alpha \in A^{r}(M),  \tag{3.7}\\
& \alpha \neg \beta=-(\alpha, \beta) \quad \text { if } \quad \alpha, \beta \in A^{1}(M). \tag{3.8}
\end{align*}

\begin{lemma} If $\alpha \in A_{b}(\mathcal{F})^{\perp} \cap A^{r, 0}(\mathcal{F})$ and $\beta \in A_{b}(\mathcal{F})$, then $\alpha \wedge \beta, \alpha \neg \beta \in A_{b}(\mathcal{F})^{\perp}$.
\end{lemma}

\begin{proof} By (1.18), (1.19) and (3.5) it is enough to prove that $\alpha \wedge \beta \in A_{b}(\mathcal{F})^{\perp}$. But for $\gamma \in A_{b}(\mathcal{F})$, by (3.6) and (3.7) we have
$$
\langle\alpha \wedge \beta, \gamma\rangle=-(-1)^{r(\tau-1) / 2}\langle\beta, \alpha \neg \gamma\rangle=0.
$$
\end{proof}

If $(T, \mathcal{H})$ is any representative of the holonomy pseudogroup of $\mathcal{F}$, then $\left(A_{b}(\mathcal{F}), d_{b}\right)$ can be canonically identified with the differential algebra $\left(A_{\mathcal{H}}(T), d_{T}\right)$ of $\mathcal{H}$-invariant differential forms on $T$ ([13]). The transverse Riemannian structure of $\mathcal{F}$ is defined by an $\mathcal{H}$-invariant Riemannian metric on $T$, and $*_Q$ on $A_{b}(\mathcal{F})$ corresponds to the star operator $*_T$ on $A_{\mathcal{H}}(T)$. Let $\delta_{T}$ be the formal adjoint of $d_{T}$ on $T$. We have
\begin{equation*}
\delta_{T}\left(A_{\mathcal{H}}(T)\right) \subset A_{\mathcal{H}}(T), \tag{3.9}
\end{equation*}
and $\delta_{T \mid A_{\mathcal{H}}(T)}$ can be extended to an operator on $A^{\cdot, 0}(M)$ (also denoted by $\delta_{T}$) by defining
\begin{equation*}
\delta_{T}=(-1)^{q r+q+1} *_Q d_{1,0} *_Q \quad \text { on } \quad A^{r, 0}(M). \tag{3.10}
\end{equation*}
We will denote $D_{T}=d_{T}+\delta_{T}$ and $\Delta_{T}=D_{T}^{2}$.

The proof of the following result is similar to that of Proposition 4.13 in [16].

\begin{proposition} On $A^{*, 0}(M)$ we have
$$
\delta_{-1,0}=\delta_{T}-\kappa \neg .
$$
\end{proposition}

By Theorem 2.1, $d_{b}$ has an adjoint operator $\delta_{b}$ given by the composition of $\delta_{-1,0}$ with the orthogonal projection onto $A_{b}(\mathcal{F})$. We will denote $D_{b}=d_{b}+\delta_{b}$ and $\Delta_{b}=D_{b}^{2}$. The following result describes $\delta_{b}$.

\begin{corollary} $\delta_{b}=\delta_{T}-\kappa_{b} \neg$.
\end{corollary}

\begin{proof}This is a direct consequence of Lemma 3.1 and Proposition 3.2.
\end{proof}

\begin{corollary}  $\quad \delta_{-1,0} \nu=*_{Q} \kappa$, and $\delta_{b} \nu=*_{Q} \kappa_{b}$.
\end{corollary}

\begin{proof}
By Proposition 3.2 and (3.5) we have
$$
\delta_{-1,0} \nu=\delta_{T} \nu-\kappa \neg \nu=*_{Q}\left(\kappa \wedge *_{Q} \nu\right)=*_{Q} \kappa,
$$
which also implies $\delta_{b} \nu=*_{Q} \kappa_{b}$ by (1.18) and (1.19).
\end{proof}

\begin{corollary}   $\quad d_{b} \kappa_{b}=0$.
\end{corollary}

\begin{proof} By Corollary 3.3 and Corollary 3.4 we have
$$
\begin{aligned}
0 & =\delta_{b}^{2} \nu=\delta_{b} *_{Q} \kappa_{b} \\
& \left.=\left(\delta_{T}-\kappa_{b}\right\urcorner\right) * Q \kappa_{b}=\delta_{T} *_{Q} \kappa_{b} \\
& =*_{Q} d_{b} \kappa_{b}.
\end{aligned}
$$
\end{proof}

Therefore $\kappa_{b}$ defines a class $\left[\kappa_{b}\right] \in H_{b}^{1}(\mathcal{F})$, even without the assumptions of orientability of $\mathcal{F}$ and $M$.

\section{Invariance of $\kappa_{b}$ by Changing the Orthogonal Complement of $T \mathcal{F}$}

In this section the following theorem will be proved.

\begin{theorem} Let $\mathcal{F}$ be a Riemannian foliation on a compact manifold $M$ with a bundle-like metric. Then the basic component of the mean curvature is invariant under changes of the bundle-like metric keeping the same transverse Riemannian structure and the same metric along the leaves.
\end{theorem}

To prove this theorem we can assume that $\mathcal{F}$ is oriented and transversally oriented.

Let $g$ and $g^{\prime}$ be two bundle-like metrics on $M$ defining the same transverse Riemannian structure and the same metric along the leaves. We may have $T \mathcal{F}^{\perp} \neq T \mathcal{F}^{\perp}$, where $\perp$ and $\perp^{\prime}$ are the corresponding relations of orthogonality. Both metrics define the same volume form $\nu$, but different characteristic forms, $\chi$ and $\chi^{\prime}$, respectively.

The metrics $g$ and $g^{\prime}$ define different bigradings of $A(M)$. Then for any $\alpha \in A(M)$ let $\alpha_{u, v}$ and $\alpha_{u, v}^{\prime}$ denote its $g$- and $g^{\prime}$-bihomogeneous component of bidegree $(u, v)$, respectively.

Let $\delta$ and $\delta^{\prime}$ be the $g$-and $g^{\prime}$-codifferential operator, respectively, and let $\delta_{i, j}$ and $\delta_{i, j}^{\prime}$ be their bihomogeneous components with respect to the corresponding bigradings of $A(M)$. Since $g$ and $g^{\prime}$ define the same metric along the leaves, by (1.9) and (1.12) we obtain
\begin{equation*}
\left(\delta_{0,-1} \alpha\right)_{u, v-1}^{\prime}=\delta_{0,-1}^{\prime}(\alpha_{u, v}^{\prime}) \tag{4.1}
\end{equation*}
for all $\alpha \in A^{u, v}(M)$.

\begin{lemma}$g$ and $g^{\prime}$ define the same scalar product on $A_{b}(M)$, and we have
\begin{equation*}
A_{b}(\mathcal{F})^{\perp}=A_{b}(\mathcal{F})^{\perp^{\prime}} \tag{4.2}
\end{equation*}
\end{lemma}

\begin{proof} Because $\chi$ and $\chi^{\prime}$ have the same restriction to the leaves, we have $\nu \wedge \chi=$ $\nu \wedge \chi^{\prime}$. Hence both metrics define the same scalar product on $A_{b}(\mathcal{F})$ since they define the same transverse Riemannian structure.

(4.2) follows easily from (4.1) and (2.2) of Theorem 2.1.
\end{proof}

Suppose that $\mathcal{F}$ is T.P. and take a $g$-unit vector field $X \in \mathcal{X}(M, \mathcal{F}) \cap \Gamma\left(T \mathcal{F}^{\perp}\right)$. Then there exists some $Y \in \mathcal{X}(\mathcal{F})$ such that $X^{\prime}=X+Y \in \Gamma\left(T \mathcal{F}^{\perp}\right)$. $X^{\prime}$ is a $g^{\prime}$-unit vector. Let $\beta \in A^{0,1}(M)$ be the $g$-dual form of $Y$, and $E_{1}, \ldots, E_{p}$ a local orthonormal frame for $T \mathcal{F}$ that is leafwise synchronous at some point $x \in M$ (see [3] and [25]). Then, by (1.8) we have at the point $x$
$$
\begin{aligned}
\delta_{0,-1} \beta & =-\sum_{i=1}^{p} \iota_{E_{i}} \stackrel{\circ}{\nabla}_{E_{i}} \beta=-\sum_{i=1}^{p} E_{i} \beta\left(E_{i}\right)=-\sum_{i=1}^{p} E_{i}\left(Y, E_{i}\right) \\
& =-\sum_{i=1}^{p}\left(\left(\nabla_{E_{i}} Y, E_{i}\right)+\left(Y, \nabla_{E_{i}} E_{i}\right)\right)=-\sum_{i=1}^{p}\big((\nabla_{E_{i}} Y, E_{i})+(Y, \stackrel{\circ}{\nabla}_{E_{i}} E_{i})\big) \\
& =-\sum_{i=1}^{p}\left(\nabla_{E_{i}} Y, E_{i}\right)=-\sum_{i=1}^{p}\left(E_{i},\left[E_{i}, Y\right]\right)=\sum_{i=1}^{p}\left(\left(\left[X^{\prime}, E_{i}\right], E_{i}\right)-\left(\left[X, E_{i}\right], E_{i}\right)\right) \\
& =\kappa^{\prime}\left(X^{\prime}\right)-\kappa(X)=\left(\kappa^{\prime}-\kappa\right)(X).
\end{aligned}
$$
Therefore, by (2.2) and (2.4) we obtain that $\kappa^{\prime}-\kappa \in A_{b}(\mathcal{F})^{\perp}$. So the theorem follows in this case by Lemma 4.2.

In the general case we consider again the $O(q)$-principal bundle of transverse orthonormal frames with the transverse Levi-Civita connection. Then, with the notation of Section 2, the characteristic form of $\hat{\mathcal{F}}$ is $\hat{\chi}=\pi^{*} \chi$. So by (3.1) we have
\begin{equation*}
\kappa=\rho_{\pi}(\hat{\kappa}), \tag{4.3}
\end{equation*}
where $\hat{\kappa}$ is the mean curvature of $\hat{\mathcal{F}}$. Then the theorem follows by (2.10) and the definition of $\rho_{\mathcal{F}}$.

We also have the following reciprocal of Theorem 4.1.

\begin{proposition} Let $\mathcal{F}$ be a Riemannian foliation on a compact manifold $M$ with a bundle-like metric $g$. For any $\gamma \in A^{1,1}(M)$ there exists a bundle-like metric $g^{\prime}$ on $M$ such that $\kappa_{b^{\prime}}^{\prime}=\kappa_{b}$ and $\kappa_{o^{\prime}}^{\prime}=\kappa_{o}+\delta_{0,-1}(\gamma)$. Moreover, $g^{\prime}$ can be chosen defining the same metric along the leaves and the same transverse Riemannian structure as $g$.
\end{proposition}

\begin{proof}As before we can assume that $\mathcal{F}$ is oriented and T.P. Then take $X_{1}, \ldots, X_{q} \in$ $\mathcal{X}(M, \mathcal{F}) \cap \Gamma\left(T \mathcal{F}^{\perp}\right)$ defining an orthonormal transverse parallelism, and let $\alpha_{1}, \ldots, \alpha_{q} \in A_{b}^{1}(\mathcal{F})$ be their dual forms.

By (2.4) there exist $\beta_{1}, \ldots, \beta_{q} \in A^{0,1}(M)$ such that
\begin{equation*}
\gamma=\alpha_{1} \wedge \beta_{1}+\ldots+\alpha_{q} \wedge \beta_{q} . \tag{4.4}
\end{equation*}
For each $i=1, \ldots, q$ let $Y_{i} \in \mathcal{X}(\mathcal{F})$ be the dual vector field of $\beta_{i}$. The vector fields $X_{i}^{\prime}=X_{i}+Y_{i}$ define a vector subbundle $Q^{\prime} \subset T M$ complementary to $T \mathcal{F}$. Consider the metric $g^{\prime}$ defined as the orthogonal sum of the restriction of $g$ to $T \mathcal{F}$ and the metric on $Q^{\prime}$ given by the transverse Riemannian structure. This change of metric is of the type we have studied, so
\begin{equation*}
\kappa^{\prime}-\kappa=\delta_{0,-1}(\gamma) \tag{4.5}
\end{equation*}
as we saw. Hence, by Lemma 4.2 the result follows.
\end{proof}

\textbf{Remark.}
The proof of Proposition 4.3 is related to Sullivan's purification of $p-$ forms ([29]) in the following way. Let $\chi$ be the $g$-characteristic form of $\mathcal{F}$. Let $\chi_{1}$ be the purification of $\chi-*_{\mathcal{F}}\left(\alpha_{1} \wedge \beta_{1}\right)$, and for each $i=2, \ldots, q$ let $\chi_{i}$ be the purification of $\chi_{i-1}-* \mathcal{F}\left(\alpha_{i} \wedge \beta_{i}\right)$. Then it is easy to check that $\chi_{q}$ is the characteristic form of $g^{\prime}$.

\section{Invariance of $\left[\kappa_{b}\right]$}
Assume that $\mathcal{F}$ is oriented and transversally oriented. Now let $g$ and $g^{\prime}$ be two bundle-like metrics on $M$ defining the same transverse Riemannian structure and such that
\begin{equation*}
T \mathcal{F}^{\perp}=T \mathcal{F}^{\perp^{\prime}}. \tag{5.1}
\end{equation*}
Then the corresponding characteristic forms satisfy the relation
\begin{equation*}
\chi=e^{f} \chi^{\prime} \tag{5.2}
\end{equation*}
for some $f \in \mathcal{C}^{\infty}(M)$. By (3.1) and (5.2) we obtain that the corresponding mean curvatures satisfy
\begin{equation*}
\kappa^{\prime}=\kappa-d_{1,0}(f). \tag{5.3}
\end{equation*}
For any $\omega \in A(M)$, let $\omega_{b^{\prime}}$ and $\omega_{o^{\prime}}$ be the $g^{\prime}$-orthogonal projections of $\omega$ into $A_{b}(\mathcal{F})$ and $A_{b}(\mathcal{F})^{\perp^{\prime}}$, respectively. From (5.3) we cannot derive $\left[\kappa_{b}\right]=\left[\kappa_{b^{\prime}}^{\prime}\right]$ diretly, for $A_{b}(\mathcal{F})^{\perp}$ may be different from $A_{b}(\mathcal{F})^{\perp^{\prime}}$. In fact, we have the following result.

\begin{proposition} Let $\delta$ and $\delta^{\prime}$ be the $g$-and $g^{\prime}$-formal adjoints of $d$, respectively.\\
We have:\\
{\rm i)} $\delta_{0,-1}^{\prime}\left(A^{\cdot, 1}(M)\right)=e^{-f} \delta_{0,-1}\left(A^{\cdot, 1}(M)\right)$;\\
{\rm ii)} $A_{b}(\mathcal{F})^{\perp^{\prime}}=e^{-f} A_{b}(\mathcal{F})^{\perp}$;\\
{\rm iii)} if $\mathcal{F}$ is basic, then $A_{b}(\mathcal{F})^{\perp^{\prime}}=A_{b}(\mathcal{F})^{\perp}$.
\end{proposition}

\begin{proof} Clearly i) implies ii) by (2.2) of Theorem 2.1. Moreover, ii) implies iii) by Lemma 3.1. So all we have to prove is i).

Let $*_{\mathcal{F}}$ and $*_{\mathcal{F}}^{\prime}$ be the $g$- and $g^{\prime}$-leafwise star operators, respectively. For $\alpha \in$ $A^{\cdot, 1}(M)$ we have
$$
\begin{aligned}
\delta_{0,-1}^{\prime}(\alpha) & =*_{\mathcal{F}}(\delta_{0,-1}^{\prime}(\alpha) \wedge \chi)=e^{-f} *_\mathcal{F}(\delta_{0,-1}^{\prime}(\alpha) \wedge \chi^{\prime}) \\
& =e^{-f} *_\mathcal{F}^{\prime}*_{\mathcal{F}} \delta_{0,-1}^{\prime}(\alpha)=-e^{-f} *_\mathcal{F} d_{0,1} *_{\mathcal{F}}^{\prime} \alpha \\
& =(-1)^{p+1} e^{-f} \delta_{0,-1} *_\mathcal{F}{*}^{\prime}_{\mathcal{F}} \alpha
\end{aligned}
$$
by (5.2) and (1.13).
\end{proof}

By ii) of Proposition 2.2 we have
\begin{equation*}
(e^{f})_{b}=(e^{f_{o}})_{b} e^{f_{b}}. \tag{5.4}
\end{equation*}
Hence, by (3.1), (5.2) and (5.4) we obtain
\begin{align*}
d_{1,0}(\chi^{\prime})&= e^{f_{b}}\left(d_{1,0}((e^{f_{o}})_{b})+(e^{f_{o}})_{b} d_{1,0}(f_{b})\right. \\
&\phantom{={}}{} \left.-(e^{f_{o}}t)_{b} \kappa_{b}\right) \wedge \chi-(e^{f_{o}})_{b} e^{f_{b}} \kappa_{o} \wedge \chi+d_{1,0}((e^{f})_{o} \chi) . \tag{5.5}
\end{align*}
On the other hand, by (3.1) and (5.2) we also have
\begin{equation*}
d_{1,0}(\chi^{\prime})=-e^{f_{b}}\left((e^{f_{o}})_{b}+(e^{f_{o}})_{o}\right) \kappa_{b^{\prime}}^{\prime} \wedge \chi-\kappa_{o^{\prime}}^{\prime} \wedge \chi^{\prime}. \tag{5.6}
\end{equation*}
But by Theorem 2.1, Lemma 3.1, and (1.13) we obtain that
\begin{equation*}
\left(e^{f_{o}}\right)_{b} e^{f_{b}} \kappa_{o} \wedge \chi, \quad d_{1,0}\left(\left(e^{f}\right)_{o} \chi\right), e^{f_{b}}\left(e^{f_{o}}\right)_{o} \kappa_{b}^{\prime}, \wedge \chi,\quad\text{and}\quad\kappa_{o^{\prime}}^{\prime} \wedge \chi^{\prime}\quad\text{belong to}\quad\overline{\mathcal{B}_{p}}. \tag{5.7}
\end{equation*}
Then, by Theorem 2.1 and (1.13),
\begin{equation*}
-(e^{f_{o}})_{b} \kappa_{b^{\prime}}^{\prime}=d_{1,0}((e^{f_{o}})_{b})+(e^{f_{o}})_{b} d(f_{b})-(e^{f_{o}})_{b} \kappa_{b}, \tag{5.8}
\end{equation*}
which implies
\begin{equation*}
\kappa_{b}-\kappa_{b^{\prime}}^{\prime}=d(f_{b}+\ln((e^{f_{o}})_{b})) \tag{5.9}
\end{equation*}
since $(e^{f_{o}})_{b}>0$ by i) of Proposition 2.2. Therefore,
\begin{equation*}
\left[\kappa_{b}\right]=\left[\kappa_{b^{\prime}}\right] \quad \text { in } \quad H_{b}^{1}(\mathcal{F}) ; \tag{5.10}
\end{equation*}
i.e. $\left[\kappa_{b}\right]$ is invariant under this type of metric changes.

Now, by (3.1) the mean curvature depends only on $T \mathcal{F}^{\perp}$ and the volume form along the leaves. So all the possibilities for the mean curvature can be obtained by using the two types of modifications of the bundle-like metric that we have seen.

Moreover, by (2.2) of Theorem 2.1, (1.9), (1.12), and (4.2) of Lemma 4.2, $A_{b}(\mathcal{F})^{\perp}$ depends only on the metric along the leaves. So we obtain all the possibilities for $A_{b}(\mathcal{F})^{\perp}$ with the last type of changes of the bundle-like metric.

Therefore, neither $\kappa$ nor $A_{b}(\mathcal{F})^{\perp}$ depend on the transverse Riemannian structure. Thus we obtain all the possibilities for $\kappa_{b}$ with the above two types of metric changes.

Finally, the assumptions about orientability can be omitted with standard arguments, obtaining the following result.

\begin{theorem} For a Riemannian foliation $\mathcal{F}$ on a compact manifold, all the bundle-like metrics define the same class $\left[\kappa_{b}\right]$ in $H_{b}^{1}(\mathcal{F})$.
\end{theorem}

So $\left[\kappa_{b}\right]$ is an invariant of the foliation that will be denoted by $\xi(\mathcal{F})$. Furthermore, by using (5.9) with $f$ basic we obtain the following result.

\begin{proposition} Every element of $\xi(\mathcal{F})$ can be achieved as the basic component of the mean curvature of some bundle-like metric defining any given transverse Riemannian structure.
\end{proposition}

\section{Tautness Characterization}

\begin{theorem} For a Riemannian foliation $\mathcal{F}$ on a compact manifold $M$ with a bundle-like metric, we have the orthogonal decomposition
\begin{equation*}
A_{b}(\mathcal{F})=\operatorname{Ker}\left(\Delta_{b}\right) \oplus \operatorname{Im}\left(d_{b}\right) \oplus \operatorname{Im}\left(\delta_{b}\right) \tag{6.1}
\end{equation*}
with respect to the restriction of the scalar product of $A(M)$.
\end{theorem}

This result can be proved with the same arguments as the basic Hodge decomposition of [18] or [23] by using $\kappa_{b}$ instead of $\kappa$, thus without assuming that $\kappa$ is basic. It has also an easier proof by checking that $-i D_{b}$ verifies the hypotheses of Lemma 2.1 of [5] in the Hilbert space completion of $A_{b}(\mathcal{F}) \otimes \mathbb{C}$, and applying directly the abstract Hodge decomposition of [3].

A basic Hodge decomposition with the same generality was also proved in [8]. But in that result the direct summands are orthogonal with respect to another scalar product induced by Molino's structure theorems. This difference will be used in the following consequence.

\begin{corollary} For a transversally oriented Riemannian foliation $\mathcal{F}$ on a compact manifold, we have $\xi(\mathcal{F})=0$ iff $H_{b}^{q}(\mathcal{F}) \neq 0$. Moreover, in this case we have $H_{b}^{q}(\mathcal{F}) \neq$ 0 iff $H_{b}^{r}(\mathcal{F}) \cong H_{b}^{q-r}(\mathcal{F})$ for each integer $r$.
\end{corollary}

\begin{proof} Clearly $\xi(\mathcal{F})=0$ implies $H_{b}^{q}(\mathcal{F}) \neq 0$ for we can choose a bundle-like metric so that $\kappa_{b}=0$, obtaining $\delta_{b} \nu=\delta_{T} \nu=0$. The reciprocal follows because if $H_{b}^{q}(\mathcal{F}) \neq$ 0 , then there exists some bundle-like metric for which $\nu$ is orthogonal to $d\left(A_{b}^{q-1}(\mathcal{F})\right)$ ([19]), and thus $\kappa_{b}=0$ by Corollary 3.4.

The duality of $H_{b}(\mathcal{F})$ when $H_{b}^{q}(\mathcal{F}) \neq 0$ was proved in [8].
\end{proof}

A foliation $\mathcal{F}$ on a manifold $M$ is said to be taut when there exists a Riemannian metric on $M$ for which all the leaves are minimal submanifolds. If $\mathcal{F}$ is Riemannian and taut, then we can always choose a bundle-like metric for which the leaves are minimal submanifolds. X.~Masa has proved in [19] that when $\mathcal{F}$ is Riemannian and oriented, and $M$ is compact and oriented, $\mathcal{F}$ is taut iff $H_{b}^{q}(\mathcal{F}) \neq 0$. The following result will be useful to weaken the above orientation hypotheses.

\begin{lemma} Let $\mathcal{F}$ be a Riemannian foliation on a compact manifold M. Let $\pi: \tilde{M} \rightarrow M$ be a finite covering of $M$. Then $\mathcal{F}$ is taut iff $\tilde{\mathcal{F}}=\pi^{*} \mathcal{F}$ is taut.
\end{lemma}

\begin{proof} Suppose that $\tilde{\mathcal{F}}$ is taut, and let $\tilde{g}$ be a minimizing metric on $\tilde{M}$. Clearly we can assume that $\tilde{\mathcal{F}}$ is oriented to prove that $\mathcal{F}$ is taut, so $\tilde{\mathcal{F}}$ has a $\tilde{g}$-characteristic form $\tilde{\chi}$.

Let $\Gamma$ be the group of deck transformations of $\pi$. We have the homomorphism $\varepsilon: \Gamma \rightarrow\{1,-1\}$ defined by $\varepsilon(\gamma)=1$ if $\gamma$ preserves the orientation of $\tilde{\mathcal{F}}$, and $\varepsilon(\gamma)=-1$ if $\gamma$ changes the orientation of $\tilde{\mathcal{F}}$. The p-form
$$
\sum_{\gamma \in \Gamma} \varepsilon(\gamma) \gamma^{*}(\tilde{\chi})
$$
restricts to a volume form along the leaves, and we can thus consider its Sullivan purification $\bar{\chi}$ ([29]). Clearly $\gamma^{*}(\bar{\chi})=\varepsilon(\gamma) \bar{\chi}$ for all $\gamma \in \Gamma$, so $\bar{\chi}$ is the $\left(\pi^{*} g\right)$-characteristic form of $\tilde{\mathcal{F}}$ for some metric $g$ on $M$.

Moreover, from (3.1) the $\tilde{g}$-minimality of the leaves of $\tilde{\mathcal{F}}$ is equivalent to $d \tilde{\chi}\left(X_{1}, \ldots\right.$, $\left.X_{p+1}\right)=0$ whenever $p$ of the vectors $X_{i}$ are tangent to the leaves (see [10]). It is easy to check that this property is also satisfied by $\bar{\chi}$, obtaining that the leaves of $\mathcal{F}$ are $g$-minimal submanifolds, and $\mathcal{F}$ is thus taut.

Reciprocally, if $\mathcal{F}$ is taut, then $\tilde{\mathcal{F}}$ is obviously taut.
\end{proof}

\textbf{Remark.} The non-trivial part of the above result cannot be proved by a direct average process of a minimizing metric on $\tilde{M}$ because the minimality of the leaves is not preserved by sums of metrics. There is no generalization of the above lemma to the case of arbitrary infinite coverings, otherwise we would obtain that every Riemannian foliation is taut by considering the universal covering of $M$, which is not true ([4]).

From Corollary 6.2, Lemma 6.3, and the tautness theorem proved by X.~Masa, we obtain the following result.

\begin{theorem} A Riemannian foliation $\mathcal{F}$ on a compact manifold is taut iff $\xi(\mathcal{F})=0$. Moreover, when $\mathcal{F}$ is transversally oriented, it is taut iff $H_{b}^{q}(\mathcal{F}) \neq 0$.
\end{theorem}

Therefore, the tautness theorem of X.~Masa was generalized to transversally oriented Riemannian foliations on compact manifolds. The orientation of $\mathcal{F}$ is not needed.

In particular, $\mathcal{F}$ is taut when $H^{1}(M, \mathbb{R})=0$, for example when $M$ is simply connected. This was already proved in [10], and also in [30] when the mean curvature is basic. (See also [6].) We also have that $\mathcal{F}$ is taut when the fundamental group of $M$ is finite.

In general, any vanishing theorem for $H_{b}^{1}(\mathcal{F})$ has a tautness corollary, as in the following example. Let $\mathcal{F}$ be a Riemannian foliation on a manifold with a complete bundle-like metric, and let $(T, \mathcal{H})$ be any representative of the holonomy pseudogroup. Let $\operatorname{Ric}_{T}$ be the Ricci curvature of the corresponding $\mathcal{H}$-invariant Riemannian metric on $T$. J. Hebda has proved in [14] that if there is a positive lower bound on $\operatorname{Ric}_{T}$, then the leaf space of $\mathcal{F}$ is compact and $H_{b}^{1}(\mathcal{F})=0$. This result was proved using relations between the focal points of the leaves and the geodesics orthogonal to the leaves. For compact manifolds it was also proved in [20] using the Weitzenböck formula and the maximum principle. Therefore we obtain the following.

\begin{corollary} Let $\mathcal{F}$ be a Riemannian foliation on a compact manifold. If $\operatorname{Ric}_{T}$ is strictly positive, then $\mathcal{F}$ is taut.\end{corollary}

For codimension one we obtain stronger consequences.

\begin{corollary} Every Riemannian foliation of codimension one on a compact manifold is taut.\end{corollary}

\begin{proof} Let $\mathcal{F}$ be such a foliation. By Lemma 6.3 we can assume that $\mathcal{F}$ is transversally oriented. Then Corollary 6.2 implies $\xi(\mathcal{F})=0$, and the result follows by Theorem 6.4.
\end{proof}

\begin{corollary} Every compact manifold $M$ with a Riemannian foliation $\mathcal{F}$ of codimension one has infinite fundamental group.
\end{corollary}

\begin{proof} Clearly we can assume that $\mathcal{F}$ is transversally oriented, thus $0 \neq H_{b}^{1}(\mathcal{F}) \subset$ $H^{1}(M, \mathbb{R})$ by Theorem 6.4.
\end{proof}

\begin{corollary} Every closed leaf of a transversally oriented Riemannian foliation of codimension one on a compact manifold is cut by a closed transversal curve.
\end{corollary}

\begin{proof}  Apply Corollary 3 of [29].
\end{proof}

\begin{corollary} Every compact semi-simple Lie group has no Lie subgroups of codimension one.
\end{corollary}

\begin{proof}  It is well-known that such a Lie group has finite fundamental group, and the left translates of any Lie subgroup define a Riemannian foliation.
\end{proof}

\end{document}


\begin{abstract}
A small correction of a lemma in the paper ``The basic component of the mean curvature of Riemannian foliations.''
\end{abstract}

\maketitle


Let $\FF$ be a Riemannian foliation on a compact manifold $M$. Let $g,g'$ be bundle-like metrics on $M$ defining the same transverse Riemannian structure and the same metric along the leaves. In \cite{Alv1992c}, from Equations~(4.1) and (2.2) of Theorem~2.1, it follows that Equation~(4.2) holds for the corresponding orthogonal complements of the basic complex $A_b(\mathcal{F})$ in the space $A^{{\cdot},0}(M)$ of differential forms that vanish at vectors tangent to the leaves. In general, Equation~(4.2) is not true on the whole of $A(M)$, as stated in Lemma~4.2. So the correct statement of Lemma~4.2  is the following.
\setcounter{section}{4}
\setcounter{thm}{1}
\setcounter{equation}{1}
\begin{lem}
$g$ and $g'$ define the same scalar product on $A_b(\mathcal{F})$, and we have
\begin{equation}
A_b(\mathcal{F})^\perp\cap A^{{\cdot},0}(M)=A_b(\mathcal{F})^{\perp'}\cap A^{{\cdot},0}(M)\;.
\end{equation}
\end{lem}

In Section~5 of the paper, in the three paragraphs preceding Theorem~5.2, in three places, $A_b(\mathcal{F})^\perp$ should be replaced with $A_b(\mathcal{F})^\perp\cap A^{{\cdot},0}(M)$. 

These changes do not affect the rest of the paper, because the mean curvature form $\kappa$ is in $A^{1,0}(M)$.

To see that two metrics $g$ and $g'$ as above may actually have $A_b(\mathcal{F})^\perp\ne A_b(\mathcal{F})^{\perp'}$, consider the following simple example.

\begin{example*}
Consider $M=S^1\times S^1=\{(x,y):x,y\in \mathbb{R}\mod 2\pi \}$, with the foliation given by the first coordinate.  In the standard product flat metric $g$, the basic one-forms are 
$\{f(y)\,dy\}$, where $f$ is periodic. The space $A_b(\mathcal{F})^\perp\cap A^1(M)$ is $\{h_1(x,y)\,dy+h_2(x,y)\,dx\}$, where $h_1$ and $h_2$ are doubly periodic and $h_1$ integrates to $0$ over $x$. We now define the new metric $g'$ by declaring that
$\{\partial_x, \partial_y+\partial_x\}$ is a global orthonormal frame, and then the inner products on the leaves and on $TM\slash T\mathcal{F}$ are the same.  Then the dual basis is 
$\{dx-dy, dy\}$ and is orthonormal for $g'$. But then $dx-dy\in A_b(\mathcal{F})^{\perp'}$, but $dx-dy\notin A_b(\mathcal{F})^{\perp}$.
\end{example*}

\bibliographystyle{amsplain}


\providecommand{\bysame}{\leavevmode\hbox to3em{\hrulefill}\thinspace}
\providecommand{\MR}{\relax\ifhmode\unskip\space\fi MR }
\providecommand{\MRhref}[2]{%
  \href{http://www.ams.org/mathscinet-getitem?mr=#1}{#2}
}
\providecommand{\href}[2]{#2}